\documentclass[dvips,noinfoline]{article}

\RequirePackage[lnms]{imsart}
\RequirePackage{imsartprelims}
\RequirePackage[colorlinks]{hyperref}

\pubyear{2008}
\volume{2}
\volumetitle{Probability and Statistics:\break Essays in Honor of\break David A. Freedman}
\firstpage{1}
\lastpage{7}
\setattribute{printed}{text}{Printed in Lithuania}

\startlocaldefs
\setattribute{journal}{name}{Collections}
\endlocaldefs

\begin{document}

\begin{titlepage}
\editor{Deborah Nolan and Terry Speed, Editors}
\end{titlepage}

\begin{copyrightpage}
\LCCN{2007939122}
\ISBN[13]{978-0-940600-74-4}
\ISBN[10]{0-940600-74-9}
\ISSN{1939-4039}
\serieseditor{Anthony Davison}
\treasurer{Rong Chen}
\executivedirector{Elyse Gustafson}

\end{copyrightpage}


\makeatletter
\gdef\doi@base{http://arXiv.org/abs/}
\makeatother

\begin{contents}[doi]
\contentsline{begintocitem}{}{}
\contentsline{jobname}{imscoll1pr}{}
\contentsline{doi}{0806.4441}{}
\contentsline{title}{Preface}{v}
\contentsline{author}{Deborah Nolan and Terry Speed}{v}
\contentsline{endtocitem}{}{}

\contentsline{begintocitem}{}{}
\contentsline{jobname}{lnms55PR}{}
\contentsline{doi}{0806.4441}{}
\contentsline{title}{Contributors}{vi}
\contentsline{author}{}{vi}
\contentsline{endtocitem}{}{}

\contentsline{section}{Probability}{}

\contentsline{begintocitem}{}{}
\contentsline{jobname}{imscoll201}{}
\contentsline{doi}{0805.2801}{}
\contentsline{title}{Probability theory and its models}{1}
\contentsline{author}{Paul Humphreys}{1}
\contentsline{endtocitem}{}{}

\contentsline{begintocitem}{}{}
\contentsline{jobname}{imscoll202}{}
\contentsline{doi}{0805.2808}{}
\contentsline{title}{Dutch book in simple multivariate normal prediction: Another look}{12}
\contentsline{author}{Morris L. Eaton}{12}
\contentsline{endtocitem}{}{}

\contentsline{begintocitem}{}{}
\contentsline{jobname}{imscoll203}{}
\contentsline{doi}{0706.2013}{}
\contentsline{title}{A transient Markov chain with finitely many cutpoints}{24}
\contentsline{author}{Nicholas James, Russell Lyons and Yuval Peres}{24}
\contentsline{endtocitem}{}{}

\contentsline{begintocitem}{}{}
\contentsline{jobname}{imscoll204}{}
\contentsline{doi}{math/0602091}{}
\contentsline{title}{Moments of convex distribution functions and completely alternating sequences}{30}
\contentsline{author}{Alexander Gnedin and Jim Pitman}{30}
\contentsline{endtocitem}{}{}

\contentsline{begintocitem}{}{}
\contentsline{jobname}{imscoll205}{}
\contentsline{doi}{math/0509270}{}
\contentsline{title}{Brownian motion on disconnected sets, basic hypergeometric functions, and some continued fractions of Ramanujan}{42}
\contentsline{author}{Shankar Bhamidi, Steven N. Evans, Ron Peled and Peter Ralph}{42}
\contentsline{endtocitem}{}{}

\contentsline{section}{Statistics}{}

\contentsline{begintocitem}{}{}
\contentsline{jobname}{imscoll206}{}
\contentsline{doi}{0805.2829}{}
\contentsline{title}{Characteristics of hand and machine-assigned scores to college students' answers to open-ended tasks}{76}
\contentsline{author}{Stephen P. Klein}{76}
\contentsline{endtocitem}{}{}

\contentsline{begintocitem}{}{}
\contentsline{jobname}{imscoll207}{}
\contentsline{doi}{0805.2835}{}
\contentsline{title}{Alternative formulas for synthetic dual system estimation in the 2000 census}{90}
\contentsline{author}{Lawrence Brown and Zhanyun Zhao}{90}
\contentsline{endtocitem}{}{}

\contentsline{begintocitem}{}{}
\contentsline{jobname}{imscoll208}{}
\contentsline{doi}{0805.2838}{}
\contentsline{title}{On the history and use of some standard statistical models}{114}
\contentsline{author}{E. L. Lehmann}{114}
\contentsline{endtocitem}{}{}

\contentsline{begintocitem}{}{}
\contentsline{jobname}{imscoll209}{}
\contentsline{doi}{0805.2840}{}
\contentsline{title}{Counting the homeless in Los Angeles County}{127}
\contentsline{author}{Richard Berk, Brian Kriegler and Donald Ylvisaker}{127}
\contentsline{endtocitem}{}{}

\contentsline{begintocitem}{}{}
\contentsline{jobname}{imscoll210}{}
\contentsline{doi}{0805.2845}{}
\contentsline{title}{Statistical adjustment for a measure of healthy lifestyle doesn't yield the truth about hormone therapy}{142}
\contentsline{author}{Diana B. Petitti and Wansu Chen}{142}
\contentsline{endtocitem}{}{}

\contentsline{begintocitem}{}{}
\contentsline{jobname}{imscoll211}{}
\contentsline{doi}{0805.3008}{}
\contentsline{title}{Multiple tests of association with biological annotation metadata}{153}
\contentsline{author}{Sandrine Dudoit, S\"{u}nd\"{u}z Kele\c {s} and Mark J. van der Laan}{155}
\contentsline{endtocitem}{}{}

\contentsline{begintocitem}{}{}
\contentsline{jobname}{imscoll212}{}
\contentsline{doi}{0805.2879}{}
\contentsline{title}{Multivariate data analysis: The French way}{219}
\contentsline{author}{Susan Holmes}{219}
\contentsline{endtocitem}{}{}

\contentsline{begintocitem}{}{}
\contentsline{jobname}{imscoll213}{}
\contentsline{doi}{0805.2881}{}
\contentsline{title}{The future of census coverage surveys}{234}
\contentsline{author}{Kenneth Wachter}{234}
\contentsline{endtocitem}{}{}

\contentsline{begintocitem}{}{}
\contentsline{jobname}{imscoll214}{}
\contentsline{doi}{0805.3019}{}
\contentsline{title}{Three months journeying of a Hawaiian monk seal}{246}
\contentsline{author}{David R. Brillinger, Brent S. Stewart and Charles L. Littnan}{246}
\contentsline{endtocitem}{}{}

\contentsline{begintocitem}{}{}
\contentsline{jobname}{imscoll215}{}
\contentsline{doi}{0805.3043}{}
\contentsline{title}{Projection pursuit for discrete data}{265}
\contentsline{author}{Persi Diaconis and Julia Salzman}{265}
\contentsline{endtocitem}{}{}

\contentsline{begintocitem}{}{}
\contentsline{jobname}{imscoll216}{}
\contentsline{doi}{0805.3027}{}
\contentsline{title}{DNA Probabilities in \textit{People v.
Prince}: When are racial and ethnic statistics relevant?}{289}
\contentsline{author}{David H. Kaye}{289}
\contentsline{endtocitem}{}{}

\contentsline{begintocitem}{}{}
\contentsline{jobname}{imscoll217}{}
\contentsline{doi}{0805.3032}{}
\contentsline{title}{Testing earthquake predictions}{302}
\contentsline{author}{Brad Luen and Philip B. Stark}{302}
\contentsline{endtocitem}{}{}

\contentsline{begintocitem}{}{}
\contentsline{jobname}{imscoll218}{}
\contentsline{doi}{0805.3034}{}
\contentsline{title}{Curse-of-dimensionality revisited: Collapse of the particle filter in very large scale systems}{316}
\contentsline{author}{Thomas Bengtsson, Peter Bickel and Bo Li}{316}
\contentsline{endtocitem}{}{}

\contentsline{begintocitem}{}{}
\contentsline{jobname}{imscoll219}{}
\contentsline{doi}{0805.3040}{}
\contentsline{title}{Higher order influence functions and minimax estimation of nonlinear\hfill\break functionals}{335}
\contentsline{author}{James Robins, Lingling Li, Eric Tchetgen and Aad van der Vaart}{335}
\contentsline{endtocitem}{}{}

\end{contents}

\begin{preface}

  \begin{frontmatter}

    \title{Preface}

  \end{frontmatter}

  \thispagestyle{plain}
This volume is our tribute to David A. Freedman, whom we regard as one
of the great statisticians of our time. He received his B.Sc. degree
from McGill University and his Ph.D. from Princeton, and joined the
Department of Statistics of the University of California, Berkeley,
in 1962, where, apart from sabbaticals, he has been ever since.

In a career of over 45 years, David has made many fine contributions
to probability and statistical theory, and to the application of
statistics. His early research was on Markov chains and martingales,
and two topics with which he has had a lifelong fascination:
exchangeability and De Finetti's theorem, and the consistency of Bayes
estimates. His asymptotic theory for the bootstrap was also highly
influential. David was elected to the American Academy of Arts and
Sciences in 1991, and in 2003 he received the John J. Carty Award for
the Advancement of Science from the U.S. National Academy of Sciences.

In addition to his purely academic research, David has extensive
experience as a consultant, including working for the Carnegie
Commission, the City of San Francisco, and the Federal Reserve, as
well as several Departments of the U.S. Government--Energy, Treasury,
Justice, and Commerce. He has testified as an expert witness on
statistics in a number of law cases, including \textit{Piva v. Xerox}
(employment discrimination), \textit{Garza v. County of Los Angeles} (voting
rights), and \textit{New York v. Department of Commerce} (census adjustment).

Lastly, he is an exceptionally good writer and teacher, and his many
books and review articles are arguably his most important contribution
to our subject. His widely used elementary text \textit{Statistics}, written
with R. Pisani and R. Purves, now in its 4th edition, is rightly
regarded as a classic introductory exposition, while his second text
\textit{Statistical Models} (2005) is set to become just as successful in
its field.

The roles of theoretical researcher, consultant, and expositor are not
disjoint aspects of David's personality, but fully integrated
ones. For over 20 years now, he has been writing extensively on
statistical modeling.  He has contributed to theory, and prepared
illuminating expositions and given penetrating critiques of old and
new models and methods in a wide range of contexts. The result is a
quite remarkable body of research on the theory and application of
statistics, particularly to the decennial U.S. census, the social
sciences (especially econometrics, political science and the law), and
epidemiology.  These themes are reflected in this volume of papers by
friends and colleagues of David's. We'd like to thank him for his
wonderful body of work, and to wish him well for the future.
\bigskip
{\flushright

\begin{tabular}{l@{}}
Deborah Nolan\\
Terry Speed\\
\end{tabular}

}

\end{preface}

\newpage
\begin{contributors}
\begin{theindex}

  \item Bengtsson, T., \textit {Bell Labs}
  \item Berk, R., \textit {University of California, Los Angeles}
  \item Bhamidi, S., \textit {University of California, Berkeley}
  \item Bickel, P., \textit {University of California, Berkeley}
  \item Brillinger, D. R., \textit {University of California, Berkeley}
  \item Brown, L., \textit {University of Pennsylvania}

  \indexspace

  \item Chen, W., \textit {Kaiser Permanente Southern California}

  \indexspace

  \item Diaconis, P., \textit {Stanford University}
  \item Dudoit, S., \textit {University of California, Berkeley}

  \indexspace

  \item Eaton, M. L., \textit {University of Minnesota, Minneapolis}
  \item Evans, S. N., \textit {University of California, Berkeley}

  \indexspace

  \item Gnedin, A., \textit {Rijksuniversiteit Utrecht}

  \indexspace

  \item Holmes, S., \textit {Stanford University}
  \item Humphreys, P., \textit {University of Virginia}

  \indexspace

  \item James, N., \textit {Berkeley, California}

  \indexspace

  \item Kaye, D. H., \textit {Arizona State University}
  \item Kele\c {s}, S., \textit {University of Wisconsin, Madison}
  \item Klein, S. P., \textit {GANSK \& Associates}
  \item Kriegler, B., \textit {University of California, Los Angeles}

  \indexspace

  \item Lehmann, E. L., \textit {University of California, Berkeley}
  \item Li, B., \textit {Tsinghua University}
  \item Li, L., \textit {Harvard School of Public Health}
  \item Littnan, C. L., \textit {Pacific Islands Fisheries Science Center}
  \item Luen, B., \textit {University of California, Berkeley}
  \item Lyons, R., \textit {Indiana University}

  \indexspace

  \item Peled, R., \textit {University of California, Berkeley}
  \item Peres, Y., \textit {University of California, Berkeley, and Microsoft Research}
  \item Petitti, D. B., \textit {University of Southern California}
  \item Pitman, J., \textit {University of California, Berkeley}

  \indexspace

  \item Ralph, P., \textit {University of California, Berkeley}
  \item Robins, J., \textit {Harvard School of Public Health}

  \indexspace

  \item Salzman, J., \textit {Stanford University}
  \item Stark, P. B., \textit {University of California, Berkeley}
  \item Stewart, B. S., \textit {Hubbs-SeaWorld Research Institute}

  \indexspace

  \item Tchetgen, E., \textit {Harvard School of Public Health}

  \indexspace

  \item van der Laan, M. J., \textit {University of California, Berkeley}
  \item van der Vaart, A., \textit {Vrije Universiteit}

  \indexspace

  \item Wachter, K., \textit {University of California, Berkeley}

  \indexspace

  \item Ylvisaker, D., \textit {University of California, Los Angeles}

  \indexspace

  \item Zhao, Z., \textit {Mathematica Policy Research, Inc.}

\end{theindex}

\end{contributors}

%

\end{document}